\documentclass{article}
\makeatletter
\newcommand{\artstat@writing}{writing}
\newcommand{\artstat@publish}{publish}
\newcommand{\artstat@check}{check}
\newcommand{\artstat}{}         
\renewcommand{\artstat}{writing}
\renewcommand{\artstat}{publish}
\newcommand{\ifwriting}[1]{\ifx\artstat\artstat@writing #1\fi}
\newcommand{\ifpublish}[1]{\ifx\artstat\artstat@publish #1\fi}
\newcommand{\ifcheck}[1]{\ifx\artstat\artstat@check #1\fi}
\makeatother

\usepackage{amssymb}
\usepackage{amsthm}
\usepackage{amsmath}
\usepackage{mathtools}
\usepackage{enumitem}
\usepackage[capitalize]{cleveref}
\usepackage{tikz}
\usepackage{pxpgfmark}          
\usetikzlibrary{calc}
\usetikzlibrary{intersections}
\usetikzlibrary{positioning}
\usepackage{tikz-cd}

\newcommand{\out}{20,10}
\newcommand{\Vout}{12,5}
\newcommand{\margin}{1,1}
\newcommand{\holeone}{10,3}
\newcommand{\Vholeone}{40:10}
\newcommand{\holetwo}{12,3}
\newcommand{\Vholetwo}{40:10}
\newcommand{\vcwidth}{3}
\newcommand{\torusscale}{0.1}

\tikzset{torus/.pic = {
    \begin{scope}[scale=\torusscale]
      \coordinate (out) at (\out);
      \coordinate (Vout) at (\Vout);
      \path[use as bounding box] let \p0 = ($(out)+(\margin)$) in (\p0) -- (-\x0,-\y0);
      \path let \p0 = ($(out)+(\margin)$) in
      coordinate (-right) at (\x0, 0)
      coordinate (-left) at (-\x0, 0);

      \draw[name path=outpath] let
      \p{out}=(out),
      \p{Vout}=(Vout)
      in (\x{out},0)
      .. controls +(0,\y{Vout}) and +(\x{Vout},0) .. (0,\y{out})
      .. controls +(-\x{Vout},0) and +(0,\y{Vout}) .. (-\x{out},0)
      .. controls +(0,-\y{Vout}) and +(-\x{Vout},0) .. (0,-\y{out})
      .. controls +(\x{Vout},0) and +(0,-\y{Vout}) .. cycle
      ;

      \begin{scope}
        \coordinate (hole1) at (\holeone);
        \coordinate (Vhole1) at (\Vholeone);
        \draw[clip,name path=cuppath] let \p{hole1}=(hole1), \p{Vhole1}=(Vhole1) in
        (\p{hole1})
        .. controls +(-\x{Vhole1},-\y{Vhole1}) and +(\x{Vhole1},-\y{Vhole1}) .. (-\x{hole1},\y{hole1});

        \coordinate (hole2) at (\holetwo);
        \coordinate (Vhole2) at (\Vholetwo);
        \draw let \p{hole2}=(hole2), \p{Vhole2}=(Vhole2) in
        (\x{hole2},-\y{hole2})
        .. controls +(-\x{Vhole2},\y{Vhole2}) and +(\x{Vhole2},\y{Vhole2}) .. (-\x{hole2},-\y{hole2});
      \end{scope}

      \path[name path=axis] let \p0=(\out) in (0,0) -- ($2*(0,-\y0)$);
      \path[name intersections={of=outpath and axis}]
        coordinate (out-bottom) at (intersection-1);
      \path[name intersections={of=cuppath and axis}]
        coordinate (cup-bottom) at (intersection-1);
      \draw (out-bottom) .. controls +(\vcwidth,0) and +(\vcwidth,0) .. (cup-bottom);
      \draw[dashed] (out-bottom) .. controls +(-\vcwidth,0) and +(-\vcwidth,0) .. (cup-bottom);
    \end{scope}
  }}
\tikzset{collapse/.pic = {
    \begin{scope}[scale=\torusscale]
      \newcommand{\collapsept}{0,-8}
      \newcommand{\Vcolin}{60:10}
      \newcommand{\Vcolout}{-30:10}
      \newcommand{\Vholeonecol}{70:10}

      \path[use as bounding box] let \p0 = ($(\out)+(\margin)$) in (\p0) -- (-\x0,-\y0);

      \path [fill] (\collapsept) circle (1);

      \draw let
      \p{out}=(\out),
      \p{Vout}=(\Vout),
      \p{hole1}=(\holeone),
      \p{Vhole1}=(\Vholeonecol),
      \p{collapsept}=(\collapsept),
      \p{Vcolout}=(\Vcolout),
      \p{Vcolin}=(\Vcolin)
      in
      (\p{hole1})                   
      .. controls +(-\x{Vhole1},-\y{Vhole1}) and +(\p{Vcolin}) ..
      (\p{collapsept})              
      .. controls +(\p{Vcolout}) and +(0,-\y{Vout}) ..
      (\x{out},0)                   
      .. controls +(0,\y{Vout}) and +(\x{Vout},0) ..
      (0,\y{out})                   
      .. controls +(-\x{Vout},0) and +(0,\y{Vout}) ..
      (-\x{out},0)                  
      .. controls +(0,-\y{Vout}) and +(-\x{Vcolout},\y{Vcolout}) ..
      (\p{collapsept})              
      .. controls +(-\x{Vcolin},\y{Vcolin}) and +(\x{Vhole1},-\y{Vhole1}) ..
      (-\x{hole1},\y{hole1})        
      ;
      \clip let
      \p{hole1}=(\holeone),
      \p{Vhole1}=(\Vholeonecol),
      \p{collapsept}=(\collapsept),
      \p{Vcolin}=(\Vcolin)
      in
      (\p{hole1})                   
      .. controls +(-\x{Vhole1},-\y{Vhole1}) and +(\p{Vcolin}) ..
      (\p{collapsept})              
      .. controls +(-\x{Vcolin},\y{Vcolin}) and +(\x{Vhole1},-\y{Vhole1}) ..
      (-\x{hole1},\y{hole1})        
      ;

      \draw let \p{hole2}=(\holetwo), \p{Vhole2}=(\Vholetwo) in
      (\x{hole2},-\y{hole2})
      .. controls +(-\x{Vhole2},\y{Vhole2}) and +(\x{Vhole2},\y{Vhole2}) .. (-\x{hole2},-\y{hole2});
    \end{scope}
  }}
\tikzset{sphere/.pic = {
    \begin{scope}[scale=\torusscale]
      \newcommand{\separateptright}{1,-8}
      \newcommand{\collapseptright}{1.5,-8}
      \newcommand{\Vsepin}{70:7}
      \newcommand{\Vsepout}{-70:6}
      \newcommand{\Vholeonesep}{80:10}
      \newcommand{\VoutAdd}{0,2}

      \path[use as bounding box] let \p0 = ($(\out)+(\margin)$) in (\p0) -- (-\x0,-\y0);

      \path [fill] let \p0=(\collapseptright) in
      (\x0,\y0) circle (1)
      (-\x0,\y0) circle (1);

      \draw let
      \p{out}=(\out),
      \p{Vout}=(\Vout),
      \p{hole1}=(\holeone),
      \p{Vhole1}=(\Vholeonesep),
      \p{separateptright}=(\separateptright),
      \p{separateptleft}=(-\x{separateptright}, \y{separateptright}),
      \p{Vsepout}=(\Vsepout),
      \p{Vsepin}=(\Vsepin),
      \p{VoutAdd}=(\VoutAdd)
      in
      (\p{hole1})                   
      .. controls +(-\x{Vhole1},-\y{Vhole1}) and +(\p{Vsepin}) ..
      (\p{separateptright})         
      .. controls +(\p{Vsepout}) and +(0,-\y{Vout}-\y{VoutAdd}) ..
      (\x{out},0)                   
      .. controls +(0,\y{Vout}) and +(\x{Vout},0) ..
      (0,\y{out})                   
      .. controls +(-\x{Vout},0) and +(0,\y{Vout}) ..
      (-\x{out},0)                  
      .. controls +(0,-\y{Vout}-\y{VoutAdd}) and +(-\x{Vsepout},\y{Vsepout}) ..
      (\p{separateptleft})          
      .. controls +(-\x{Vsepin},\y{Vsepin}) and +(\x{Vhole1},-\y{Vhole1}) ..
      (-\x{hole1},\y{hole1})        
      ;
      \clip let
      \p{hole1}=(\holeone),
      \p{Vhole1}=(\Vholeonesep),
      \p{separateptright}=(\separateptright),
      \p{separateptleft}=(-\x{separateptright}, \y{separateptright}),
      \p{Vsepin}=(\Vsepin)
      in
      (\p{hole1})                   
      .. controls +(-\x{Vhole1},-\y{Vhole1}) and +(\p{Vsepin}) ..
      (\p{separateptright})              
      --
      (\p{separateptleft})
      .. controls +(-\x{Vsepin},\y{Vsepin}) and +(\x{Vhole1},-\y{Vhole1}) ..
      (-\x{hole1},\y{hole1})        
      ;

      \draw let \p{hole2}=(\holetwo), \p{Vhole2}=(\Vholetwo) in
      (\x{hole2},-\y{hole2})
      .. controls +(-\x{Vhole2},\y{Vhole2}) and +(\x{Vhole2},\y{Vhole2}) .. (-\x{hole2},-\y{hole2});

    \end{scope}
  }}

\newtheoremstyle{mydescription}
  {}
  {}
  {}
  {}
  {}
  {}
  { }
  {(\thmnumber{#2})}

\numberwithin{equation}{section}

\theoremstyle{definition}
\newtheorem{definition}[equation]{Definition}
\theoremstyle{plain}

\newtheorem{proposition}[equation]{Proposition}
\newtheorem{theorem}[equation]{Theorem}
\newtheorem{corollary}[equation]{Corollary}

\newtheorem{problem}[equation]{Problem}
\theoremstyle{remark}

\theoremstyle{mydescription}
\newtheorem{mydescription}[equation]{MyDescription}

\crefname{equation}{}{}
\crefname{enumi}{}{}
\crefname{mydescription}{}{}
\crefformat{mydescription}{(#2#1#3)}

\usepackage{todonotes}
\presetkeys{todonotes}{color=orange!60}{}

\makeatletter
\tikzstyle{inlinenotestyle} = [
    draw=\@todonotes@currentbordercolor,
    fill=\@todonotes@currentbackgroundcolor,
    line width=0.5pt,
    inner sep = 0.8 ex,
    rounded corners=4pt]

\renewcommand{\@todonotes@drawInlineNote}{%
        {\begin{tikzpicture}[remember picture,baseline=(current bounding box.base)]%
            \draw node[inlinenotestyle,font=\@todonotes@sizecommand, anchor=base,baseline]{%
              \if@todonotes@authorgiven%
                {\noindent \@todonotes@sizecommand \@todonotes@author:\,\@todonotes@text}%
              \else%
                {\noindent \@todonotes@sizecommand \@todonotes@text}%
              \fi};%
           \end{tikzpicture}}}%
\newcommand{\mytodo}[1]{\@todo[inline]{#1}}%
\makeatother

\makeatletter
\ifwriting{%
  \usepackage[color]{showkeys}
  \SK@def\cref#1{\SK@\SK@@ref{#1}\SK@cref{#1}} 
  \definecolor{refkey}{rgb}{0.7, 0.8, 0.5}
  \definecolor{labelkey}{rgb}{0, 0.7, 0.5}
}
\makeatother
\ifcheck{\presetkeys{todonotes}{disable}{}}
\ifpublish{\presetkeys{todonotes}{disable}{}}

\newcommand{\zw}[1]{\makebox[0pt][l]{#1}}

\newcommand{\id}{{\rm id}}

\newcommand{\res}{\mathrm{res}}

\newcommand{\K}{\mathbb K}
\newcommand{\Q}{\mathbb Q}
\newcommand{\Z}{\mathbb Z}

\newcommand{\cochain}[2][*]{C^{#1}(#2)}
\newcommand{\cohom}[2][*]{H^{#1}(#2)}
\newcommand{\homol}[2][*]{H_{#1}(#2)}

\newcommand{\im}{\operatorname{Im}}

\renewcommand{\hom}{\operatorname{Hom}}
\newcommand{\ext}{\operatorname{Ext}}
\newcommand{\tor}{\operatorname{Tor}}
\renewcommand{\lim}{{\varprojlim}}

\newcommand{\pt}{\mathrm{pt}}

\newcommand{\piratq}{\bigoplus_n\pi_n(M)\otimes\Q}
\newcommand{\pirat}{\bigoplus_n\pi_n(M)\otimes\K}
\newcommand{\pirats}{\pi_*(M)\otimes\K}
\newcommand{\evenpart}[1]{#1^{\rm even}}
\newcommand{\oddpart}[1]{#1^{\rm odd}}

\newcommand{\shriek}[1]{{#1}_!}
\newcommand{\incl}{\mathrm{incl}}
\newcommand{\comp}{\mathrm{comp}}
\newcommand{\diag}{\Delta}
\newcommand{\inclconst}{c}

\newcommand{\tpow}[1]{^{\otimes #1}}
\newcommand{\susp}[2][]{s^{#1}#2}

\newcommand{\selfwedge}[1]{{#1}_\vee}
\newcommand{\selfconn}[1]{{#1}_\#}

\newcommand{\lprod}{\mu}
\newcommand{\lcop}{\delta}
\newcommand{\bprod}{\mu}
\newcommand{\bcop}{\delta}
\newcommand{\dprod}{\mu^\vee}
\newcommand{\dcop}{\delta^\vee}
\DeclareMathOperator{\map}{Map}

\newcommand{\spheresp}[2][]{S^{#1}#2}

\title{Nontrivial example of the composition of the brane product and coproduct on Gorenstein spaces}
\author{Shun Wakatsuki}
\date{}

\begin{document}
\maketitle
\begin{abstract}
  We give an example of a space with the nontrivial composition of
  the brane product and the brane coproduct,
  which we introduced in a previous article.
  \todo{composition followed by?}
\end{abstract}

\makeatletter
\ifx\artstat\artstat@writing
  \listoftodos
\fi
\makeatother

\section{Introduction}
Chas and Sullivan \cite{chas-sullivan} introduced the loop product
$\lprod\colon\homol{LM\times LM}\to \homol[*-m]{LM}$ on the homology of the free loop space
$LM = \map(S^1, M)$ of a connected closed oriented manifold $M$ of dimension $m$.
Constructing a 2-dimensional topological quantum field theory without counit,
Cohen and Godin \cite{cohen-godin} generalized this product to other string operations,
including the loop coproduct
$\lcop\colon\homol{LM}\to\homol[*-m]{LM\times LM}$.
But Tamanoi \cite{tamanoi} showed that
any string operation corresponding to a positive genus surface is trivial.
In particular, the composition $\lprod \circ \lcop$ is trivial.
There are many attempts
to find nontrivial and interesting operations.

F\'elix and Thomas \cite{felix-thomas09} generalized the loop product and coproduct
to the case $M$ is a Gorenstein space.
A Gorenstein space is a generalization of a manifold
in the point of view of Poincar\'e duality.
For example,
connected closed oriented manifolds,
classifying spaces of connected Lie groups,
and their Borel constructions \todo{Lie acting on mfd} are Gorenstein spaces.
Moreover, any 1-connected space $M$ with $\piratq$ of finite dimension is a Gorenstein space.
In spite of this huge generalization,
string operations remain to tend to be trivial.
For example, the loop product $\lprod$ is trivial over a field of characteristic zero
for the classifying space of a connected Lie group \cite[Theorem 14]{felix-thomas09}.

\begin{problem}[\cite{felix-thomas09}]\todo{引用}
  \label{problem:involutive}
  Is there a Gorenstein space
  such that the composition $\lprod\circ\lcop$ is nontrivial?
\end{problem}
This is the Gorenstein counterpart of the above result due to Tamanoi.
But such an example is not found.

Sullivan and Voronov \cite[Part I, Chapter 5]{cohen-hess-voronov}
generalized the loop product to the sphere space $\spheresp[k]{M}=\map(S^k, M)$ for $k\geq 1$.
This product is called the brane product.

The brane coproduct, a generalization of the loop coproduct to the sphere spaces,
is constructed by the author \cite{wakatsuki18:toappear}
in the case where the rational homotopy group $\piratq$ is of finite dimension.
In the construction,
we assume the ``finiteness'' of the dimension of the $(k-1)$-fold based loop space $\Omega^{k-1} M$
as a Gorenstein space.
Moreover,
the product and the coproduct were generalized
to the mapping spaces from manifolds, by means of connected sums.

\newcommand{\src}{S}
\newcommand{\inclorig}{p}
\newcommand{\comporig}{q}
Here we briefly review the brane product and coproduct.
See \cref{section:braneOperations} for details.
Let $\K$ be a field of characteristic zero,
$\src$ an oriented manifold of dimension $k$ with two disjoint base points, and
$M$ a $k$-connected $m$-dimensional $\K$-Gorenstein space with $\pirat$ of finite dimension.
Denote the ``connected sum'' and ``wedge sum'' of $\src$ with itself
along the two base points
by $\selfconn\src$ and $\selfwedge\src$, respectively.
Note that, by the definition of the connected sum,
we have the canonical inclusion
$S^{k-1}\hookrightarrow \selfconn\src$
and the quotient map
$\comporig\colon \selfconn\src \to (\selfconn\src)/S^{k-1} = \selfwedge\src$.
Similarly we have
$S^0=\pt\coprod\pt \hookrightarrow \src$
and
$\inclorig\colon \src \to S/S^0 = \selfwedge\src$.
Hence we have the following diagram

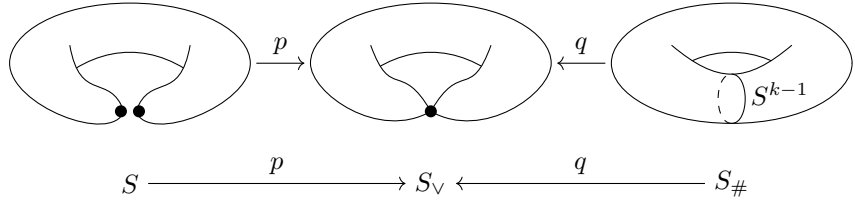
\begin{tikzpicture}[scale=0.8]
  \path (0,1.6);                
  \path pic[transform shape] at (0,0) {sphere};
  \path pic[transform shape] at (5,0) {collapse};
  \path pic[transform shape] at (9.99999,0) {torus}; 
  \draw[->] (2.1,0) -- (2.9,0) node[midway,above] {$\inclorig$};
  \draw[->] (7.9,0) -- (7.1,0) node[midway,above] {$\comporig$};
  \node (S) at (0,-2) {$\src$};
  \node (wedge) at (5,-2) {$\selfwedge\src$};
  \node (conn) at (10,-2) {$\selfconn\src$};
  %
  \draw[->] (S.east) -- (wedge.west) node[midway,above] {$\inclorig$};
  \draw[->] (conn.west) -- (wedge.east) node[midway,above] {$\comporig$};
  \node at (10.8,-0.5) {$S^{k-1}$};
\end{tikzpicture}

\noindent
and its dual

\begin{equation}
  \label{equation:inclAndComp}
  M^\src \xleftarrow{\incl}
  M^{\selfwedge\src} \xrightarrow{\comp}
  M^{\selfconn\src},
\end{equation}
where the maps $\incl$ and $\comp$ are induced by $\inclorig$ and $\comporig$.

Using this diagram,
we can construct two operations, $\src$-brane product $\bprod_{\src}$ and coproduct $\bcop_{\src}$:
\begin{align*}
  \bprod_{\src}\colon& \homol{M^\src} \rightarrow \homol[*-m]{M^{\selfconn\src}} \\
  \bcop_{\src}\colon& \homol{M^{\selfconn\src}} \rightarrow \homol[*-\bar{m}]{M^\src}.
\end{align*}

Note that,
if $T$ and $U$ are oriented $k$-manifolds
and we take $\src=T\coprod U$ with one base point on $T$ and the other on $U$,
then $\bprod_{\src}$ and $\bcop_{\src}$ have the form
\begin{align*}
  \bprod_{T\coprod U}\colon& \homol{M^T\times M^U} \rightarrow \homol[*-m]{M^{T\#U}} \\
  \bcop_{T\coprod U}\colon& \homol{M^{T\#U}} \rightarrow \homol[*-\bar{m}]{M^T\times M^U}.
\end{align*}
Moreover, if we take $T=U=S^1$,
then $\bprod_{S^1\coprod S^1}$ and $\bcop_{S^1\coprod S^1}$ coincide with the usual loop product and coproduct, respectively.
Hence the $\src$-brane product and coproduct are generalizations of the loop product and coproduct.


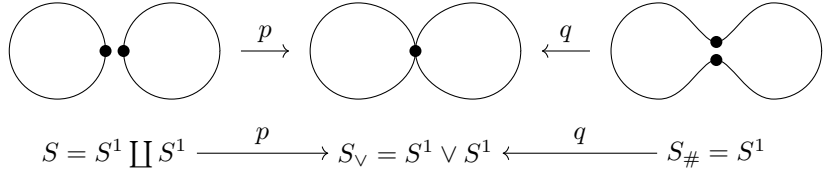
\newcommand{\radius}{8}
\newcommand{\centerx}{9.5}

\makeatletter
\@ifundefined{circlerate}{\newcommand{\circlerate}{0.55228}}{}
\makeatother

\tikzset{directsum/.pic = {
    \begin{scope}[scale=\torusscale]
      \path[draw]
      let
        \p0 = (\centerx,0),
        \n1 = {\radius},
        \n2 = {\circlerate*\n1} in
      ($(\p0)+(\n1,0)$)
        .. controls +(0,\n2) and +(\n2,0) ..
      ($(\p0)+(0,\n1)$)
        .. controls +(-\n2,0) and +(0,\n2) ..
      ($(\p0)+(-\n1,0)$)
        .. controls +(0,-\n2) and +(-\n2,0) ..
      ($(\p0)+(0,-\n1)$)
        .. controls +(\n2,0) and +(0,-\n2) ..
      cycle;
      \path[draw]
      let
        \p0 = (-\centerx,0),
        \n1 = {\radius},
        \n2 = {\circlerate*\n1} in
      ($(\p0)+(\n1,0)$)
        .. controls +(0,\n2) and +(\n2,0) ..
      ($(\p0)+(0,\n1)$)
        .. controls +(-\n2,0) and +(0,\n2) ..
      ($(\p0)+(-\n1,0)$)
        .. controls +(0,-\n2) and +(-\n2,0) ..
      ($(\p0)+(0,-\n1)$)
        .. controls +(\n2,0) and +(0,-\n2) ..
      cycle;
      \path [fill]
        ($(\centerx,0)+(-\radius,0)$) circle (1)
        ($(-\centerx,0)+(\radius,0)$) circle (1);
    \end{scope}
  }}

\tikzset{wedgesum/.pic = {
    \newcommand{\collapsept}{0,0}
    \begin{scope}[scale=\torusscale]
      \path[draw]
      let
        \p0 = (\centerx,0),
        \n1 = {\radius},
        \n2 = {\circlerate*\n1} in
      ($(\p0)+(\n1,0)$)
        .. controls +(0,\n2) and +(\n2,0) ..
      ($(\p0)+(0,\n1)$)
        .. controls +(-\n2,0) and +(0,\n2) ..
      (\collapsept)
        .. controls +(0,-\n2) and +(-\n2,0) ..
      ($(\p0)+(0,-\n1)$)
        .. controls +(\n2,0) and +(0,-\n2) ..
      cycle;
      \path[draw]
      let
        \p0 = (-\centerx,0),
        \n1 = {\radius},
        \n2 = {\circlerate*\n1} in
      (\collapsept)
        .. controls +(0,\n2) and +(\n2,0) ..
      ($(\p0)+(0,\n1)$)
        .. controls +(-\n2,0) and +(0,\n2) ..
      ($(\p0)+(-\n1,0)$)
        .. controls +(0,-\n2) and +(-\n2,0) ..
      ($(\p0)+(0,-\n1)$)
        .. controls +(\n2,0) and +(0,-\n2) ..
      cycle;
      \path [fill]
        (\collapsept) circle (1);
    \end{scope}
  }}

\tikzset{connsum/.pic = {
    \newcommand{\collapsept}{0,0}
    \newcommand{\zure}{0,1.5}
    \newcommand{\hoge}{2.5}
    \begin{scope}[scale=\torusscale]
      \path[draw]
      let
        \p0 = (\centerx,0),
        \n1 = {\radius},
        \n2 = {\circlerate*\n1},
        \p3 = (-\centerx,0),
        \n4 = {\radius},
        \n5 = {\circlerate*\n4} in
      ($(\p0)+(\n1,0)$)
        .. controls +(0,\n2) and +(\n2,0) ..
      ($(\p0)+(0,\n1)$)
        .. controls +(-\n2,0) and +(\hoge,0) ..
      ($(\collapsept)+(\zure)$)
        .. controls +(-\hoge,0) and +(\n5,0) ..
      ($(\p3)+(0,\n4)$)
        .. controls +(-\n5,0) and +(0,\n5) ..
      ($(\p3)+(-\n4,0)$)
        .. controls +(0,-\n5) and +(-\n5,0) ..
      ($(\p3)+(0,-\n4)$)
        .. controls +(\n5,0) and +(-\hoge,0) ..
      ($(\collapsept)-(\zure)$)
        .. controls +(\hoge,0) and +(-\n2,0) ..
      ($(\p0)+(0,-\n1)$)
        .. controls +(\n2,0) and +(0,-\n2) ..
      cycle;
      \path [fill]
        ($(\collapsept)+(\zure)$) circle (1)
        ($(\collapsept)-(\zure)$) circle (1);
    \end{scope}
  }}


\begin{tikzpicture}[scale=0.8]
  \path (0,1.6);                
  \path pic[transform shape] at (0,0) {directsum};
  \path pic[transform shape] at (5,0) {wedgesum};
  \path pic[transform shape] at (9.99999,0) {connsum}; 
  \draw[->] (2.1,0) -- (2.9,0) node[midway,above] {$\inclorig$};
  \draw[->] (7.9,0) -- (7.1,0) node[midway,above] {$\comporig$};
  \newcommand{\yofdiagram}{-1.7}
  \node (S) at (0,\yofdiagram) {$\src=S^1\coprod S^1$};
  \node (wedge) at (5,\yofdiagram) {$\selfwedge\src=S^1\vee S^1$};
  \node (conn) at (10,\yofdiagram) {$\selfconn\src=S^1$};
  \draw[->] (S.east) -- (wedge.west) node[midway,above] {$\inclorig$};
  \draw[->] (conn.west) -- (wedge.east) node[midway,above] {$\comporig$};
\end{tikzpicture}

In this article,
we give examples that
the composition $\bprod \circ \bcop$ of the brane product and the brane coproduct
is nontrivial.

\newcommand{\ds}{k}
\newcommand{\dt}{{k-1}}
\begin{theorem}
  \label{theorem:compNontriv}
  Let $k$ be a positive even integer.
  Consider the case $S=S^\ds$ (and hence $\selfconn{S} = S^\dt\times S^1$).
  Let $M$ be the Eilenberg-MacLane space $K(\Z, 2n)$ with $n > k/2$.
  Then the composition $\bprod_{S^k}\circ\bcop_{S^k}$ of the $S^k$-brane product
  \begin{equation*}
    \bprod_{S^k}\colon\homol{\map(S^\ds, M)}\to\homol[*+2n-1]{\map(S^\dt\times S^1, M)}
  \end{equation*}
  and the $S^k$-brane coproduct
  \begin{equation*}
    \bcop_{S^k}\colon\homol{\map(S^\dt\times S^1, M)}\to\homol[*-2n+k-1]{\map(S^\ds, M)}
  \end{equation*}
  is nontrivial.
\end{theorem}

This gives an answer to \cref{problem:involutive}
in the context of brane operations.
Here it should be remarked that,
the composition $\bprod_{S^k}\circ\bcop_{S^k}$ corresponds to a cobordism without ``genus''.
In fact,
if we take $k=1$,
the composition $\bprod_{S^1}\circ\bcop_{S^1}$
is equal to the composition $\lcop\circ\lprod$, not $\lprod\circ\lcop$,
of the loop product $\lprod$ and coproduct $\lcop$.
\todo{genus???}

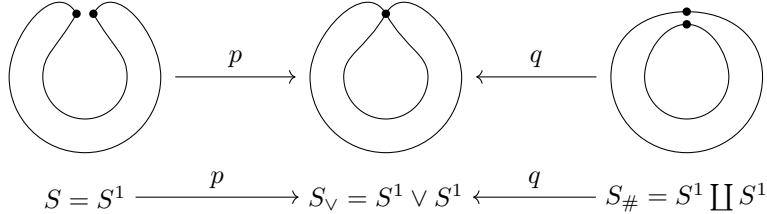
\newcommand{\outsize}{18}
\newcommand{\insize}{10}
\newcommand{\torusonescale}{0.07}

\makeatletter
\@ifundefined{circlerate}{\newcommand{\circlerate}{0.55228}}{}
\makeatother

\tikzset{torus_one/.pic = {
    \newcommand{\zure}{0,2.5}
    \begin{scope}[scale=\torusonescale]
      \path[draw]
      let
        \n0 = {\outsize},
        \n1 = {\circlerate*\n0} in
      (\n0,0)
        .. controls +(0,\n1) and +(\n1,0) ..
      ($(0,\n0)-(\zure)$)
        .. controls +(-\n1,0) and +(0,\n1) ..
      (-\n0,0)
        .. controls +(0,-\n1) and +(-\n1,0) ..
      (0,-\n0)
        .. controls +(\n1,0) and +(0,-\n1) ..
      cycle;
      \path[draw]
      let
        \n0 = {\insize},
        \n1 = {\circlerate*\n0} in
      (\n0,0)
        .. controls +(0,\n1) and +(\n1,0) ..
      ($(0,\n0)+(\zure)$)
        .. controls +(-\n1,0) and +(0,\n1) ..
      (-\n0,0)
        .. controls +(0,-\n1) and +(-\n1,0) ..
      (0,-\n0)
        .. controls +(\n1,0) and +(0,-\n1) ..
      cycle;
      \path [fill]
        ($(0,\outsize)-(\zure)$) circle (1)
        ($(0,\insize)+(\zure)$) circle (1);
    \end{scope}
  }}

\tikzset{collapse_one/.pic = {
    \newcommand{\collapsept}{0,15}
    \newcommand{\Vcolout}{60:10}
    \newcommand{\Vcolin}{-60:10}
    \begin{scope}[scale=\torusonescale]
      \path[draw]
      let
        \n0 = {\outsize},
        \n1 = {\circlerate*\n0},
        \p2 = (\Vcolout) in
      (\n0,0)
        .. controls +(0,\n1) and +(\x2,\y2) ..
      (\collapsept)
        .. controls +(-\x2,\y2) and +(0,\n1) ..
      (-\n0,0)
        .. controls +(0,-\n1) and +(-\n1,0) ..
      (0,-\n0)
        .. controls +(\n1,0) and +(0,-\n1) ..
      cycle;
      \path[draw]
      let
        \n0 = {\insize},
        \n1 = {\circlerate*\n0},
        \p2 = (\Vcolin) in
      (\n0,0)
        .. controls +(0,\n1) and +(\x2,\y2) ..
      (\collapsept)
        .. controls +(-\x2,\y2) and +(0,\n1) ..
      (-\n0,0)
        .. controls +(0,-\n1) and +(-\n1,0) ..
      (0,-\n0)
        .. controls +(\n1,0) and +(0,-\n1) ..
      cycle;
    \path [fill]
      (\collapsept) circle (1);
    \end{scope}
  }}

\tikzset{sphere_one/.pic = {
    \newcommand{\collapsept}{0,15}
    \newcommand{\zure}{2,0}
    \newcommand{\Vcolout}{60:10}
    \newcommand{\Vcolin}{-60:10}
    \begin{scope}[scale=\torusonescale]
      \path[draw]
      let
        \n0 = {\outsize},
        \n1 = {\circlerate*\n0},
        \p2 = (\Vcolout),
        \n3 = {\insize},
        \n4 = {\circlerate*\n3},
        \p5 = (\Vcolin) in
      (\n0,0)
        .. controls +(0,\n1) and +(\x2,\y2) ..
      ($(\collapsept)+(\zure)$)
        .. controls +(\x5,\y5) and +(0,\n4) ..
      (\n3,0)
        .. controls +(0,-\n4) and +(\n4,0) ..
      (0,-\n3)
        .. controls +(-\n4,0) and +(0,-\n4) ..
      (-\n3,0)
        .. controls +(0,\n4) and +(-\x5,\y5) ..
      ($(\collapsept)-(\zure)$)
        .. controls +(-\x2,\y2) and +(0,\n1) ..
      (-\n0,0)
        .. controls +(0,-\n1) and +(-\n1,0) ..
      (0,-\n0)
        .. controls +(\n1,0) and +(0,-\n1) ..
      cycle;
      \path [fill]
        ($(\collapsept)+(\zure)$) circle (1)
        ($(\collapsept)-(\zure)$) circle (1);
    \end{scope}
  }}

\begin{tikzpicture}[scale=0.8]
  \path (0,1.6);                
  \path pic[transform shape] at (0,0) {sphere_one};
  \path pic[transform shape] at (5,0) {collapse_one};
  \path pic[transform shape] at (9.99999,0) {torus_one}; 
  \draw[->] (1.5,0) -- (3.5,0) node[midway,above] {$\inclorig$};
  \draw[->] (8.5,0) -- (6.5,0) node[midway,above] {$\comporig$};
  \newcommand{\yofdiagram}{-2}
  \node (S) at (0,\yofdiagram) {$S=S^1$};
  \node (wedge) at (5,\yofdiagram) {$\selfwedge\src=S^1\vee S^1$};
  \node (conn) at (10,\yofdiagram) {$\selfconn\src=S^1\coprod S^1$};
  \draw[->] (S.east) -- (wedge.west) node[midway,above] {$\inclorig$};
  \draw[->] (conn.west) -- (wedge.east) node[midway,above] {$\comporig$};
\end{tikzpicture}

On the other hand,
the $\src$-brane coproduct is trivial in some cases.

\begin{theorem}
  \label{theorem:coprodTriv}
  Let $k$ be a positive even integer,
  and $M$ a $k$-connected (Gorenstein) space with $\pirat$ of finite dimension.
  Assume that the minimal Sullivan model of $M$ is pure
  and has at least one generator of odd degree.
  Then the $S^k$-brane coproduct is trivial for $M$.
\end{theorem}

For a connected Lie group $G$ and its closed connected subgroup $H$,
the homogeneous space $M = G / H$ satisfies the assumption
if the canonical map $\pi_*(H)\otimes\K \to \pi_*(G)\otimes\K$ is \textit{not} surjective.

By \cref{theorem:compNontriv} and \cref{theorem:coprodTriv},
we have the following corollary.

\begin{corollary}
  Let $k$ be a positive even integer,
  and $M$ a $k$-connected (Gorenstein) space with $\pirat$ of finite dimension.
  Assume that the minimal Sullivan model of $M$ is pure.
  Then the composition $\bprod_{S^k}\circ\bcop_{S^k}$ is nontrivial
  if and only if
  $M$ is a finite product $\prod K(\Z,2n_i)$ of Eilenberg-MacLane spaces of even degrees.
\end{corollary}

\cref{section:braneOperations} contains
brief background material on brane operations.
In \cref{section:modelsOfBraneOperations},
we construct rational models of the $S^k$-brane product and coproduct,
which gives a method of computation.
Next we review explicit constructions of the shriek maps
in \cref{section:explicitShriek},
which is necessary to accomplish the computation by the above models.
Finally, in \cref{section:proofOfTheorem}, we prove \cref{theorem:compNontriv} and \cref{theorem:coprodTriv}
using the above models.
\todo{要チェック}

\tableofcontents

\section{Brane operations for the mapping space from manifolds}
\label{section:braneOperations}

In this section, we review the constructions of the $\src$-brane product and coproduct
from \cite{wakatsuki18:toappear}.
Since the cochain models work well for fibrations, \todo{good for?}
we define the duals of the $\src$-brane product and coproduct at first,
and then we define the $\src$-brane product and coproduct as the duals of them.

Let $\K$ be a field of characteristic zero.
This assumption enables us to make full use of rational homotopy theory.
For the basic definitions and theorems on homological algebra and rational homotopy theory,
we refer the reader to \cite{felix-halperin-thomas01}.

\begin{definition}
  [{\cite{felix-halperin-thomas88}}]
  Let $m\in\Z$ be an integer.
  \begin{enumerate}
    \item An augmented dga (differential graded algebra) $(A,d)$ is called
      a ($\K$-\nolinebreak[4])\nolinebreak[0]\textit{Gorenstein algebra} of dimension $m$ if \todo{閉じ括弧の直前で改行}
      \begin{equation*}
        \dim \ext_A^l(\K, A) =
        \begin{cases}
          1 & \mbox{ (if $l = m$)} \\
          0 & \mbox{ (otherwise),}
        \end{cases}
      \end{equation*}
      where the field $\K$ and the dga $(A,d)$ are $(A,d)$-modules via the augmentation map and the identity map, respectively.
    \item A path-connected topological space $M$ is called a ($\K$-)\textit{Gorenstein space} of dimension $m$
      if the singular cochain algebra $\cochain{M}$ of $M$ is a Gorenstein algebra of dimension $m$.
  \end{enumerate}
\end{definition}

Here, $\ext_A(L, N)$ is defined using a semifree resolution of $(L,d)$ over $(A,d)$,
for a dga $(A,d)$ and $(A,d)$-modules $(L,d)$ and $(N,d)$.
$\tor_A(L,N)$ is defined similarly.
See \cite[Section 1]{felix-halperin-thomas01} for details of semifree resolutions.

An important example of a Gorenstein space is given by the following \lcnamecref{proposition:FinDimImplyGorenstein}.

\begin{proposition}
  [{\cite[Proposition 3.4]{felix-halperin-thomas88}}]
  \label{proposition:FinDimImplyGorenstein}
  A 1-connected topological space $M$ is a $\K$-Gorenstein space if $\pirat$ is finite dimensional.
  Similarly, a Sullivan algebra $(\wedge V, d)$ is a Gorenstein algebra if $V$ is finite dimensional.
\end{proposition}

Note that this \lcnamecref{proposition:FinDimImplyGorenstein} is proved
only for $\Q$-Gorenstein \textit{spaces} in \cite{felix-halperin-thomas88},
but the proof can be applied for any field $\K$ of characteristic zero and Sullivan \textit{algebras} over $\K$.
\todo{spaces と algebras を強調する？}

We use the following theorem \todo{$k=1$のときのFHT01に言及？引用に[FHT01 for $k=1$, Wak18 for $k\geq 2$]と書く？}
to construct the brane operations.

\begin{theorem}
  [{\cite[Theorem 12]{felix-thomas09} for $k=1$, \cite[Corollary 3.2]{wakatsuki18:toappear} for $k\geq 2$}]
  \label{theorem:extSphereSpace}
  Let $M$ be a $(k-1)$-connected (and 1-connected)
  space with $\pirat$ of finite dimension,
  for $k\geq 1$.
  Then we have an isomorphism
  \begin{equation*}
    \ext^*_{\cochain{\spheresp[k-1]{M}}}(\cochain{M}, \cochain{\spheresp[k-1]{M}}) \cong \cohom[*-\bar{m}]{M},
  \end{equation*}
  where $\bar{m}$ is the dimension of $\Omega^{k-1}M$ as a Gorenstein space.
\end{theorem}


Now we can define the $\src$-brane coproduct as follows.
Let $\src$ be an oriented manifold with two distinct base points,
$M$ a $k$-connected $m$-dimensional $\K$-Gorenstein space with $\pirat$ of finite dimension.
Consider the diagram, extending \cref{equation:inclAndComp},
\begin{equation}
  \label{equation:SCopDiagram}
  \begin{tikzcd}[row sep=large]
    M^{\selfconn\src}\arrow[d, "\res"'] & M^{\selfwedge\src} \arrow[l,"\comp"']\arrow[d]\ar[r,"\incl"] & M^\src\\
    \spheresp[k-1]{M} & M\zw. \arrow[l,"\inclconst"']
  \end{tikzcd}
\end{equation}
Here,
the square is a pullback diagram,
the map $\res$ is the restriction map to $S^{k-1}$,
and $\inclconst$ is the embedding as the constant maps.
By \cref{theorem:extSphereSpace}, we have
$\ext^{\bar{m}}_{\cochain{\spheresp[k-1]{M}}}(\cochain{M}, \cochain{\spheresp[k-1]{M}}) \cong \cohom[0]{M} \cong\K$,
hence the generator
\begin{equation*}
  \shriek\inclconst \in \ext^{\bar{m}}_{\cochain{\spheresp[k-1]{M}}}(\cochain{M}, \cochain{\spheresp[k-1]{M}})
\end{equation*}
is well-defined up to the multiplication by a non-zero scalar.
Using the map $\shriek\inclconst$ and the diagram \cref{equation:SCopDiagram},
we can define the shriek map $\shriek\comp$ as the composition
\begin{equation*}
  \begin{array}{l}
    \cohom{M^{\selfwedge\src}}
    \xleftarrow[\cong]{\mathrm{EM}} \tor^*_{\cochain{\spheresp[k-1]{M}}}(\cochain{M}, \cochain{M^{\selfconn\src}})\\
    \xrightarrow{\tor_\id(\shriek\inclconst, \id)}
    \tor^{*+\bar{m}}_{\cochain{\spheresp[k-1]{M}}}(\cochain{\spheresp[k-1]{M}}, \cochain{M^{\selfconn\src}})
    \xrightarrow[\cong]{} \cohom[*+\bar{m}]{M^{\selfconn\src}},
  \end{array}
\end{equation*}
where the map $\mathrm{EM}$ is the Eilenberg-Moore map,
which is an isomorphism since $\spheresp[k-1]{M}$ is 1-connected
(see \cite[Theorem 7.5]{felix-halperin-thomas01} for details).
By this, we define the dual of the $\src$-brane coproduct as the composition
\begin{equation*}
  \dcop_\src\colon
  \cohom[*]{M^\src}
  \xrightarrow{\incl^*} \cohom[*+\bar{m}]{M^{\selfwedge\src}}
  \xrightarrow{\shriek{\comp}} \cohom[*+\bar{m}]{M^{\selfconn\src}}.
\end{equation*}

Similarly we can define the $\src$-brane product using the generator
\begin{equation*}
  \shriek\diag \in \ext_{\cochain{M^2}}^m(\cochain{M}, \cochain{M^2})
\end{equation*}
and the diagram
\begin{equation}
  \label{equation:SProdDiagram}
  \begin{tikzcd}[row sep=large]
    M^\src \arrow[d] & M^{\selfwedge\src} \ar[l,"\incl"']\ar[d]\arrow[r,"\comp"] & M^{\selfconn\src}\\
    M\times M & M\zw. \ar[l,"\diag"']
  \end{tikzcd}
\end{equation}
Note that,
for the brane product and the \textit{loop} coproduct,
we can replace the assumption $\pirat$ is of finite dimension
with the assumption $\pirats$ is of finite type
by using \cite[Theorem 12]{felix-thomas09}
instead of \cref{theorem:extSphereSpace}.

\section{Models of the brane operations}
\label{section:modelsOfBraneOperations}
\newcommand{\torus}[1][k]{T^{(#1)}}
\newcommand{\collapse}[1][k]{U^{(#1)}}
\newcommand{\model}[1]{{\mathcal M}(#1)}
\newcommand{\m}{{\mathcal M}}

In this section,
we consider the case $S = S^k$ and
give rational models of the $S^k$-brane operations,
for an integer $k\geq 1$.
In \cref{section:proofOfTheorem},
we will prove \cref{theorem:compNontriv} and \cref{theorem:coprodTriv}
using these models.

Naito \cite{naito13} constructed a rational model of the duals of the loop product and coproduct
in terms of Sullivan models
using the torsion functor description of \cite{kuribayashi-menichi-naito}.
The author \cite{wakatsuki18:toappear} constructed a rational model of the duals of the brane product and coproduct
as a generalization of it.
Here we give a rational model of the $S^k$-brane operations
by a similar method.

\subsection{Models of spaces}

Let $M$ be a $k$-connected space with $\pirat$ of finite dimension.
Take a Sullivan model $(\wedge V, d)$ of $M$ with $V^{\leq k}=0$ and $\dim V < \infty$.
For simplicity, we sometimes denote $(\wedge V, d)$ by $\m$.
Denote $\selfconn{(S^k)} = S^{k-1}\times S^1$ by $\torus$
and $\selfwedge{(S^k)} = (S^{k-1}\times S^1) / S^{k-1}$ by $\collapse$.
For an integer $l \in \Z$, let $\susp[l]V$ be a graded module defined by $(\susp[l]V)^n=V^{n+l}$
and $\susp[l]v$ denotes the element in $\susp[l]V$ corresponding to an element $v\in V$.
Here we recall models of mapping spaces from the interval, sphere, and disk.
\newcommand{\susplow}{\susp[(k-1)]}
\newcommand{\susphigh}{\susp[(k)]}
\newcommand{\diffsphere}{d}
\newcommand{\diffdisk}{d}
\newcommand{\modelconst}{\varepsilon}
\newcommand{\diskresol}{\tilde{\varepsilon}}
\newcommand{\pathresol}{\bar{\varepsilon}}
\newcommand{\shriekinclconstrep}{\gamma}
\newcommand{\shriekdiagrep}{\eta}
\begin{mydescription}
  \label{mydescription:pathModel}
  Consider $s$ as an derivation on the algebra
  $\wedge V\tpow2\otimes \wedge \susp{V}$
  with $s \circ s = 0$.
  Define a derivation $d$ on the algebra by
  \begin{equation*}
    d(\susp v)=1\otimes v - v\otimes 1 - \sum_{i=1}^\infty\frac{(sd)^i}{i!}(v\otimes 1),
  \end{equation*}
  inductively.
  Denote the dga $(\wedge V\tpow2\otimes\wedge \susp{V}, d)$ by $\model{I}$.
  This is a Sullivan model of the path space $M^I$ ($\simeq M$).
  Moreover, define a map
  $\pathresol\colon\model{I}\to\m$
  by $\pathresol(v\otimes 1) = \pathresol(1\otimes v) = v$ and $\pathresol(\susp{v})=0$
  for $v \in V$.
  Then it is a relative Sullivan model (resolution) of the product map $\wedge V\tpow2 \wedge \susp{V}$.
  See \cite[Section 15 (c)]{felix-halperin-thomas01} or \cite[Appendix A]{wakatsuki16}
  for details.
\end{mydescription}
\begin{mydescription}
  \label{mydescription:sphereModel}
  Assume $k \geq 2$.
  Define derivations $\susplow$ and $d$ on the graded algebra $\wedge V\otimes \wedge\susp[k-1]V$ by
  \begin{align*}
    &\susplow(v)=\susp[k-1]v,\quad \susplow(\susp[k-1]v)=0, \\
    &\diffsphere(v)=dv,\ \mbox{and}\quad\diffsphere(\susp[k-1]v)=(-1)^{k-1}\susplow dv.
  \end{align*}
  Denote the dga $\wedge V\otimes \susp[k-1]V$ by $\model{S^{k-1}}$.
  This is a Sullivan model of the space $M^{S^{k-1}}$.
  See \cite[Section 5]{wakatsuki18:toappear} for details.
\end{mydescription}
\begin{mydescription}
  \label{mydescription:diskModel}
  Assume $k \geq 2$.
  Define derivations $\susphigh$ and $\diffdisk$
  on the graded algebra $\wedge V\otimes \wedge\susp[k-1]V \otimes\wedge\susp[k]V$ by
  \begin{align*}
    &\susphigh(v)=\susp[k]v,\quad\susphigh(\susp[k-1]v)=\susphigh(\susp[k]v)=0,\quad%
      \diffdisk(v)=dv, \\
    &\diffdisk(\susp[k-1]v)=\diffsphere(\susp[k-1]v),\ %
      \mbox{and}\quad\diffdisk(\susp[k]v)=\susp[k-1]v+(-1)^k\susphigh dv.
  \end{align*}
  Denote the dga $\wedge V\otimes \wedge\susp[k-1]V \otimes\wedge\susp[k]V$ by $\model{D^k}$.
  This is a Sullivan model of the space $M^{D^k}$ ($\simeq M$).
  Moreover, define a map $\diskresol\colon \model{D^k}\to \m$
  by $\diskresol(v)=v$, $\diskresol(\susp[k-1]v)=\diskresol[k]v=0$ for $v\in V$.
  Then it is a relative Sullivan model (resolution)
  of the map $\modelconst\colon\model{S^{k-1}}\to \m$,
  where $\modelconst(v)=v$ and $\modelconst(\susp[k-1]v)=0$.
  In particular, $\diskresol$ is a quasi-isomorphism.
  See \cite[Section 5]{wakatsuki18:toappear} for details.
\end{mydescription}

Next we construct models of mapping spaces
which appear in the definition of brane operations,
using the above models.

\begin{mydescription}
  \label{mydescription:torusModel}
  Since $M^{\torus} = (M^{S^{k-1}})^{S^1}$,
  we have a Sullivan model
  $\model{\torus} = (\wedge V\otimes \wedge\susp[k-1]V \otimes \wedge \susp V\otimes \wedge\susp\susp[k-1]V, d)$
  of $M^{\torus}$
  iterating the construction in \cref{mydescription:sphereModel}.
\end{mydescription}
\begin{mydescription}
  \label{mydescription:collapseModel}
  Since $\collapse$ is homotopy equivalent to $S^{k}\vee S^1$,
  the mapping space $M^{\collapse}$ is homotopy equivalent to $M^{S^k}\times_MM^{S^1}$,
  and hence we have a Sullivan model
  $\model{\collapse} = (\wedge V\otimes \wedge\susp[k]V, d) \otimes (\wedge V\otimes \wedge\susp V, d)$.
\end{mydescription}

\subsection{Models of operations}

\newcommand{\inclcochain}{\incl^*}
\newcommand{\compcochain}{\comp^*}

Here we give a model of the $S^k$-brane product and coproduct
in a similar way to \cite{naito13} and \cite{wakatsuki18:toappear}.

First we give a model of the $S^k$-brane coproduct.
Recall that the dual $\dcop_{S^k}$ of the $S^k$-brane coproduct
is the composition
\begin{equation*}
  \dcop_{S^k}\colon
  \cohom[*]{M^{S^k}}
  \xrightarrow{\incl^*} \cohom[*+\bar{m}]{M^{\collapse}}
  \xrightarrow{\shriek{\comp}} \cohom[*+\bar{m}]{M^{\torus}}.
\end{equation*}
First the map
$\incl^*\colon\cohom[*]{M^{S^k}} \to \cohom[*+\bar{m}]{M^{\collapse}}$
is induced by the canonical inclusion
$\model{S^k} \to \model{\collapse}$,
which we also denote by $\inclcochain$.
Next the map
$\shriek{\comp}\colon \cohom[*+\bar{m}]{M^{\collapse}} \to \cohom[*+\bar{m}]{M^{\torus}}$
is computed as follows.
Let
\begin{equation*}
  \shriekinclconstrep \in \hom_{\model{S^{k-1}}}(\model{D^k}, \model{S^{k-1}})
\end{equation*}
be a representative of the nontrivial element (see \cref{theorem:extSphereSpace})
\begin{align*}
  \shriek\inclconst
  \in&\ \ext^{\bar{m}}_{\cochain{\spheresp[k-1]{M}}}(\cochain{M}, \cochain{\spheresp[k-1]{M}})\\
  \cong&\ \cohom[\bar{m}]{\hom_{\model{S^{k-1}}}(\model{D^k}, \model{S^{k-1}})}.
\end{align*}
Then the map
\begin{equation*}
  \tor_\id(\shriek\inclconst, \id)\colon
  \tor^*_{\cochain{\spheresp[k-1]{M}}}(\cochain{M}, \cochain{M^{\selfconn\src}})
  \to
  \tor^{*+\bar{m}}_{\cochain{\spheresp[k-1]{M}}}(\cochain{\spheresp[k-1]{M}}, \cochain{M^{\selfconn\src}})
\end{equation*}
is induced by the cochain map
\begin{equation*}
  \shriekinclconstrep\otimes\id\colon
  \model{D^k} \otimes_{\model{{S^{k-1}}}}\model{\torus}
  \to \model{S^{k-1}} \otimes_{\model{{S^{k-1}}}}\model{\torus},
\end{equation*}
since $\model{D^k}$ is a resolution of $\m$ over $\model{S^{k-1}}$.
The map $\shriek\comp$ is computed by this combined with
the quasi-isomorphism
\begin{equation}
  \label{equation:wrongWayQuasiIsom}
  \diskresol\otimes\id\colon
  \model{D^k} \otimes_{\model{{S^{k-1}}}}\model{\torus}
  \xrightarrow[\simeq]{}
  \m \otimes_{\model{{S^{k-1}}}}\model{\torus}.
\end{equation}
Hence the dual of the $S^k$-brane coproduct is induced by the composition
\begin{equation}
  \label{equation:braneCopModel}
  \begin{array}{l}
    \model{S^k}
    \xrightarrow{\incl^*} \model{\collapse}
    \xrightarrow{\cong} \m \otimes_{\model{{S^{k-1}}}}\model{\torus} \\
    \xleftarrow[\simeq]{\diskresol\otimes\id} \model{D^k} \otimes_{\model{{S^{k-1}}}}\model{\torus} \\
    \xrightarrow{\shriekinclconstrep\otimes\id} \model{S^{k-1}} \otimes_{\model{{S^{k-1}}}}\model{\torus}
    \xrightarrow{\cong}\model{\torus}.
  \end{array}
\end{equation}

Similarly, the dual of the $S^k$-brane product is induced by the composition
\begin{equation}
  \label{equation:braneProdModel}
  \begin{array}{l}
    \model{\torus}
    \xrightarrow{\compcochain} \model{\collapse}
    \xrightarrow{\cong} \m \otimes_{\m\tpow2}(\model{I}\otimes_\m\model{S^k}) \\
    \xleftarrow[\simeq]{\pathresol\otimes\id} \model{I} \otimes_{\m\tpow2}(\model{I}\otimes_\m\model{S^k})
    \xrightarrow{\shriekdiagrep\otimes\id} \m\tpow2\otimes_{\m\tpow2}(\model{I}\otimes_\m\model{S^k}) \\
    \xrightarrow{\cong} \model{I}\otimes_\m\model{S^k}
    \xrightarrow{\pathresol\otimes\id} \m\otimes_\m\model{S^k}
    \xrightarrow{\cong} \model{S^k}
  \end{array}
\end{equation}
Here
$\shriekdiagrep \in \hom_{\m\tpow2}(\model{I}, \m\tpow2)$
is a representative of the nontrivial element
$\shriek\diag \in \ext_{\cochain{M^2}}^m(\cochain{M}, \cochain{M^2})$
and
$\compcochain\colon \model{\torus} \to \model{\collapse}$
is the canonical quotient map.

\section{Explicit construction of shriek maps}
\label{section:explicitShriek}

Models of $S^k$-brane operations are constructed in \cref{section:modelsOfBraneOperations}
using the representatives of the shriek maps $\shriekinclconstrep$ and $\shriekdiagrep$.
They are constructed by \cref{theorem:extSphereSpace},
which only states the existence of the shriek maps.
In this section,
we recall methods to construct shriek maps explicitly
from \cite{naito13}, \cite{wakatsuki16} and \cite{wakatsuki18:toappear}.

Recall the definition of a pure Sullivan algebra.
Here we denote $\evenpart{V}=\bigoplus_nV^{2n}$ and $\oddpart{V}=\bigoplus_nV^{2n+1}$.
\begin{definition}
  [c.f. {\cite[Section 32]{felix-halperin-thomas01}}]
  \label{definition:pureSullivanAlgebra}
  A Sullivan algebra $(\wedge V, d)$ with $\dim V < \infty$ is called {\it pure}
  if $d(\evenpart{V})=0$ and $d(\oddpart{V}) \subset \wedge \evenpart{V}$.
\end{definition}
In the rest of this section,
let $(\wedge V, d)$ be a pure minimal Sullivan algebra,
$\{x_1,\ldots, x_p\}$ a basis of $\evenpart{V}$, and
$\{y_1,\ldots, y_q\}$ a basis of $\oddpart{V}$.

\subsection{Construction of $\shriek\diag$}
Here we recall the description of $\shriek\diag$ in \cite{wakatsuki16},
which is a generalization of that of Naito \cite{naito13}.
Note that,
although the description holds if the Sullivan model $(\wedge V, d)$ is semi-pure
(see \cite[Definition 1.5]{wakatsuki16} for the definition),
we only refer and use it in the case $(\wedge V, d)$ is pure.

\begin{proposition}[{\cite[Theorem 5.6 (2)]{wakatsuki16}}]
  Take
  $(\wedge V\tpow2\otimes\wedge\susp V, d) = \model{I}$
  as in \cref{mydescription:pathModel}.
  If a cocycle
  $\shriekdiagrep \in \hom_{\wedge V\tpow2}(\wedge V\tpow2\otimes\wedge\susp V, \wedge V\tpow2)$
  satisfies
  \begin{equation*}
    \shriekdiagrep(\susp x_1\cdots \susp x_p)
    = (1\otimes y_1 - y_1\otimes 1)\cdots(1\otimes y_q - y_q\otimes 1),
  \end{equation*}
  then we have
  \begin{align*}
    [\shriekdiagrep] \neq 0
    \in&\ \cohom{\hom_{\wedge V\tpow2}(\wedge V\tpow2\otimes\wedge\susp V, \wedge V\tpow2)}\\
    \cong&\ \ext_{\wedge V\tpow2}(\wedge V, \wedge V\tpow2).
  \end{align*}
\end{proposition}

This proposition gives a construction of the map $\shriek\diag$.

\subsection{Construction of $\shriek\inclconst$}
Next we recall the description of $\shriek\inclconst$ in \cite{wakatsuki18:toappear}.
The following proposition gives it completely when $k$ is even.
\begin{proposition}[{\cite[Proposition 6.2]{wakatsuki18:toappear}}]
  \label{proposition:constructShriekInclconst}
  Assume that 
  $k$ is even.
  Define an element
  \begin{equation*}
    \shriekinclconstrep \in \hom_{\model{S^{k-1}}}(\model{D^k}, \model{S^{k-1}})
  \end{equation*}
  by $\shriekinclconstrep(\susp[k]y_1\cdots\susp[k]y_q)=\susp[k-1] x_1\cdots\susp[k-1] x_p$
  and $\shriekinclconstrep(\susp[k]y_{j_1}\cdots\susp[k]y_{j_l})=0$ for $l<q$.
  Then $\shriekinclconstrep$ defines a non-trivial element in
  $\ext_{\model{S^{k-1}}}(\m, \model{S^{k-1}})$.
  \todo{$\m$だったり$(\wedge V, d)$だったりするの統一しよう}
\end{proposition}

Note that, although the proposition is proved only when $k=2$ in \cite{wakatsuki18:toappear},
the same proof also applies when $k > 2$ as long as $k$ is even.

\section{Proof of \cref{theorem:compNontriv} and \cref{theorem:coprodTriv}}
\label{section:proofOfTheorem}

In this section,
we give a proof of \cref{theorem:compNontriv} and \cref{theorem:coprodTriv}
using the models constructed above.

\newcommand{\sectcop}{\varphi}
\begin{proof}[Proof of \cref{theorem:compNontriv}]
  We compute the $S^k$-brane coproduct using \cref{equation:braneCopModel}.
  Since $M=K(\Z,2n)$,
  we take the Sullivan model $(\wedge V, d) = (\wedge x, 0)$
  where $x$ is the generator of degree $2n$.
  Note that, in this case, the differentials in $\model{S^k}$ and $\model{\torus}$ are zero,
  and hence they are identified with the cohomology groups $\cohom{M^{S^k}}$ and $\cohom{M^{\torus}}$.

  By \cref{proposition:constructShriekInclconst},
  we have a representative $\shriekinclconstrep$ of the shriek map $\shriek\inclconst$ defined by
  $\shriekinclconstrep(1)=\susp[k-1]x$ and $\shriekinclconstrep((\susp[k]x)^l)=0$ for $l\geq 1$.

  Since any Sullivan algebra satisfies the lifting property for a surjective quasi-isomorphism,
  there is a section $\sectcop$ of $\diskresol\otimes\id$ in \cref{equation:wrongWayQuasiIsom},
  which is also a quasi-isomorphism.
  It is given explicitly by
  $\sectcop(1\otimes x)=1\otimes x$,
  $\sectcop(1\otimes \susp[k]x)=1\otimes \susp\susp[k-1]x$, and
  $\sectcop(1\otimes \susp x)=1\otimes \susp x$.

  Using these maps,
  we compute the composition \cref{equation:braneCopModel}. \todo{ここをちゃんと書く．5行くらい？}
  Since all maps in the composition are $\wedge V$-linear,
  it is enough to compute the image for the elements $(\susp[k]x)^n$ for $n\geq 0$.
  Applying $\incl^*$ and the section $\sectcop$ to the element,
  we have that it is mapped to
  $1\otimes(\susp\susp[k-1]x)^n \in \model{D^k} \otimes_{\model{{S^{k-1}}}}\model{\torus}$.
  Then the map $\shriekinclconstrep\otimes \id$
  send it to
  $\susp[k-1]x\otimes(\susp\susp[k-1]x)^n \in \model{S^{k-1}} \otimes_{\model{{S^{k-1}}}}\model{\torus}$.
  Hence the $S^k$-brane coproduct $\dcop_{S^k}$ is the map determined by
  $\dcop_{S^k}(\alpha) = \susp[k-1]x\iota(\alpha)$,
  where $\iota\colon \model{S^k}\to\model{\torus}$ is the algebra map defined by
  $\iota(x)=x$ and $\iota(\susp[k]x)=\susp\susp[k-1]x$.

  \newcommand{\sectprod}{\psi}
  Similarly we can compute the $S^k$-brane product.
  Define a $\wedge V\tpow2$-linear map $\shriekdiagrep\colon\model{I}\to\wedge V\tpow2$ by
  $\shriekdiagrep(1)=0$ and $\shriekdiagrep(\susp x)=1$.
  By \cref{proposition:constructShriekInclconst},
  $\shriekdiagrep$ is a representative of the shriek map $\shriek\diag$.
  We have a section $\sectprod$ of $\pathresol\otimes\id$ in \cref{equation:braneProdModel},
  which is defined by
  $\sectprod(x\otimes 1)=1\otimes (x_1\otimes 1)$,
  $\sectprod(1\otimes \susp x)=1\otimes (\susp x\otimes 1) - \susp x\otimes 1$, and
  $\sectprod(1\otimes \susp[k]x)=1\otimes (1\otimes \susp[k]x)$.
  Here we denote the element $x\otimes 1\in \model{I}$ by $x_1$.

  As a result, \todo{こっちも書く}
  the $S^k$-brane product $\dprod_{S^k}$ is the map determined by
  $\dprod_{S^k}(\beta) = 0$,
  $\dprod_{S^k}(\susp x\cdot\beta) = -\rho(\beta)$, and
  $\dprod_{S^k}(\susp[k-1]x\cdot\beta) = \dprod_{S^k}(\susp x\cdot\susp[k-1]x\cdot\beta)= 0$,
  for $\beta \in \wedge x \otimes\wedge\susp\susp[k-1]x$.
  Here $\rho\colon \wedge x\otimes\wedge\susp\susp[k-1]x\to\model{S^k}$ is the algebra map defined by
  $\rho(x)=x$ and $\rho(\susp\susp[k-1]x)=\susp[k]x$.

  Composing these two,
  we have $\bcop_{S^k}\circ\bprod_{S^k}\neq 0$.
  In fact,
  $\bcop_{S^k}\circ\bprod_{S^k}(\susp x) = -\susp[k-1]x \neq 0 \in \model{\torus}\cong\cohom{M^{\torus}}$.
  This proves the theorem.
\end{proof}

Next we prove \cref{theorem:coprodTriv}.

\begin{proof}[Proof of \cref{theorem:coprodTriv}]
  Let $(\wedge V, d)$ be the minimal Sullivan model of $M$,
  $\{x_1,\ldots, x_p\}$ a basis of $\evenpart{V}$, and
  $\{y_1,\ldots, y_q\}$ a basis of $\oddpart{V}$.
  Consider the part
  \begin{align*}
    \m \otimes_{\model{{S^{k-1}}}}\model{\torus}
    & \xleftarrow[\simeq]{\diskresol\otimes\id} \model{D^k} \otimes_{\model{{S^{k-1}}}}\model{\torus} \\
    & \xrightarrow{\shriekinclconstrep\otimes\id} \model{S^{k-1}} \otimes_{\model{{S^{k-1}}}}\model{\torus}
  \end{align*}
  in \cref{equation:braneCopModel}.
  Define a section $\sectcop$ of $\diskresol\otimes\id$ by
  $\sectcop(1\otimes v)=1\otimes v$,
  $\sectcop(1\otimes \susp v)=1\otimes \susp v$, for $v \in V$,
  $\sectcop(1\otimes\susp\susp[k-1]x_i)=1\otimes\susp\susp[k-1]x_i$, and
  $\sectcop(1\otimes\susp\susp[k-1]y_j)=1\otimes\susp\susp[k-1]y_j + (-1)^k\susp\sigma (dy_j\otimes 1)$.
  Here, in the last term $\susp\sigma (dy_j\otimes 1)$,
  $\sigma$ is the derivation which sends $v\otimes 1$ to $\susp[k]v\otimes 1$, for $v \in V$,
  and the other generators to $0$.
  The map $\susp$ is also the derivation which sends $v$ to $\susp v$, $\susp[k-1]v$ to $\susp\susp[k-1]v$,
  and others to $0$.
  \todo{check!}
  Then we have
  $\im\sectcop \subset {\mathcal N}\otimes_{\model{S^{k-1}}}\model{\torus}$,
  where
  ${\mathcal N} =\wedge V\otimes\wedge\susp[k-1]V\otimes\wedge\susp[k]\{x_1,\ldots,x_p\} \subset \model{D^k}$.
  Let $\shriekinclconstrep$ be the representative of $\shriek\inclconst$
  given by \cref{proposition:constructShriekInclconst}.
  Since $V$ has at least one generator of odd degree,
  $\shriekinclconstrep$ is zero on ${\mathcal N}$.
  This implies that the composition $(\shriekinclconstrep\otimes 1)\circ \sectcop$ is zero,
  and hence the brane coproduct $\bcop_{S^k}$ is zero.
\end{proof}

\section*{Acknowledgment}
I would like to express my gratitude to Katsuhiko Kuribayashi and Takahito Naito for productive discussions and valuable suggestions.
Furthermore, I would like to thank my supervisor Nariya Kawazumi for the enormous support and comments.
This work was supported by JSPS KAKENHI Grant Number 16J06349 and the Program for Leading Graduate School, MEXT, Japan.



\end{document}